\documentclass{amsart}
\usepackage{amsmath}
\usepackage{amssymb}
\usepackage{amsthm}
\usepackage{graphicx}
\usepackage{float}
\usepackage{pgfplots}
\usepackage{tikz}
\usetikzlibrary{arrows.meta}
\usepackage{subcaption}
\usepackage{multirow}

\theoremstyle{thmstyleone}%
\newtheorem{theorem}{Theorem}[section]%
\newtheorem{lemma}{Lemma}[section]%

\theoremstyle{thmstyletwo}%
\newtheorem{remark}{Remark}[section]%

\theoremstyle{thmstylethree}%
\newtheorem{definition}{Definition}[section]%

\theoremstyle{thmstylefour}%
\newtheorem{corollary}{Corollary}[section]%

\numberwithin{equation}{section}

\begin{document}

\title[Nonlinear eigenproblems on  Riemannian manifolds]{Nonlinear eigenvalue problems for a class of quasilinear operator on complete Riemannian manifolds}

%    author one information
\author[B. Shen]{Bin Shen}
\address{School of Mathematics, Southeast University, Nanjing 211189, P. R. China}
%\curraddr{ }
\email{shenbin@seu.edu.cn}
%\thanks{The first author is supported partially by the NNSFC (No. 12001099, 12271093)}

%    author two information
\author[Y. Zhu]{Yuhan Zhu}
\address{School of Mathematics, Southeast University, Nanjing 211189, P. R. China}
%\curraddr{ }
\email{yuhanzhu@seu.edu.cn}

%    \subjclass is required.

%    Abstract is required.
\begin{abstract}
	In this manuscript, we study the nonlinear eigenvalue problem on complete Riemannian manifolds with Ricci curvature bounded from below, to find the unknowns $\lambda$ and $u$, such that
	$$
	Qu + \lambda f(u) = 0
	$$
	where $\lambda$ is an eigenvalue of $u$, with respect to the quasilinear operator $Qu = \operatorname{div} (\mathcal{F}(u^2, |\nabla u|^2)\nabla u)$ and nonlinar function $f(\cdot)\neq 0$. We generalize the Cheng--Yau gradient estimate in \cite{shen2025feasibilitynashmoseriterationchengyautype} and demonstrate that under certain conditions, a non-zero eigenvalue gives rise to unbounded eigenfunction $u$. Our new result also covers more quasilinear equations like $p$-porous medium equation (\textit{i.e.} $\Delta_p u^q = \lambda u^r$), and generally, $\Delta_{p}\left(\sum_{i=1}^{m}a_iu^{q_i}\right)+\lambda u^r = 0$.
\end{abstract}

\maketitle
\section{Introduction}
Let $M$ be an $n$-dimensional complete Riemannian manifold and $\Omega$ is a bounded domain in $M$. The nonlinear eigenvalue problem is to find the unknowns eigenvalue $\lambda$ and eigenfunction  $u$ such that 
 \begin{align}
 \left\{\begin{aligned}
 \Delta u + \lambda f(u) = 0 &\text{ ~~~in } \Omega,\\
 u=0 &\text{ ~~~in } \partial \Omega.
 \end{aligned}\right.
 \end{align}
As shown in \cite{evans2010partial}, this problem is equivalent to find a minimizer of energy function
$$E(u) = \frac{1}{2}\int_{\Omega}|\nabla u|^2$$
under the constraint of 
$\left\{ u\in H^1_0(\Omega):\int_{\Omega}F(u) = 0\right\}$, where $F$ is a given smooth function and $f = F'$. In particular, 
there is a significant gradient estimate by Cheng and Yau \cite{cheng1975differential} that a positive harmonic function $u$ (\textit{i.e.} $\lambda = 0$) on an $n$-dimensional complete Riemannian manifold with Ricci curvature $\operatorname{Ric}\geq-(n-1)K$ for some $K\geq0$, satisfies
\begin{equation}\label{CY_esti}
\sup_{B(o,R/2)}\frac{|\nabla u|}{u}\leq c_n\frac{1+\sqrt{K}R}{R},
\end{equation}
on the geodesic ball $B(o,R)$, where $c_n$ is a constant depending only on $n$. This type of gradient estimate is vital in geometric analysis because the bound is intrinsic and invariant under different choice of local coordinates.

An interesting question is how to generalize 
the linear operator $\Delta$ to a quasilinear one
\begin{align}\label{general_nonlinear_2}
Qu := \operatorname{div}\left(\mathcal{F}(u^2,|\nabla u|^2)\nabla u\right),
\end{align}
where $\mathcal{F}: [0,\infty) \times [0,\infty)\to \mathbb{R}$. There have been many types of gradient estimates for some specific operators, such as $p$-Laplacian  $\Delta_pu = \operatorname{div}\left(|\nabla u|^{p-2}\nabla u\right)$ by B. Kotschwar and L. Ni \cite{kotschwar2009local}, and exponential Laplacian $\Delta_eu =\operatorname{div}\left(\exp\left(|\nabla u|^{2} \right)\nabla u\right)$ by J. Wu, \textit{et al.} \cite{wu2014gradient}. Also for $\varphi$-Laplacian $\Delta_{\varphi}u = \operatorname{div}\left(\varphi\left(|\nabla u|^{2} \right)\nabla u\right)$, the authors discussed several cases with different $\varphi$ and $f$ in \cite{shen2025feasibilitynashmoseriterationchengyautype}. Moreover, when $\mathcal{F} = mu^{m-1}$, the parabolic type is also related to porous medium equation ($m>1$) and fast diffusion equation ($0<m<1$):
$$\frac{\partial u}{\partial t} = \Delta u^m$$
which was discussed in \cite{huang2025gradient} by S. Huang and the first author.

In Finsler metric measure spaces, the Laplacian operator $\Delta^{\nabla u}u$ is indeed quasilinear. Thus, our generalization is also meaningful for understanding the anisotropic phenomenon in Finsler geometry. We refer the readers to \cite{xia2014local,xia2019Sharp,xia2022local} for some research in gradient estimates and eigenvalue in Finsler manifolds. Nevertheless, not all gradient estimates obtained in previous research are of the Cheng--Yau type, and there seems to be little literature directly researching on general quasilinear operators. Our goal in this series of research is to unify the above results and to cover as many equations as possible.

In the manuscript, we focus on the following nonlinear eigenproblem
\begin{align}\label{eigenproblem}
Qu = \operatorname{div} (\mathcal{F}(u^2, |\nabla u|^2)\nabla u) = - \lambda\psi(u^2)u,
\end{align}
where $\psi$ is a positive smooth function on $(0,+\infty)$. If we further assume that $\mathcal{F}$ is separable, that is, $\mathcal{F}(s,t) = a(s)\varphi(t)$, where both $a$ and $\varphi$ are $C^\infty$-functions on $[0,+\infty)$ and positive on $(0,+\infty)$,
then we have $$Qu = a(u^2)\operatorname{div} (\varphi(|\nabla u|^2)\nabla u)+ 2ua'(u^2)\varphi(|\nabla u|^2)|\nabla u|^2.$$
To further simplify, we expand the definition in \cite{shen2025feasibilitynashmoseriterationchengyautype}.
\begin{definition}\label{degree}
	For a $C^1$-function $a(t)$ on $[0,\infty)$, {the} \textbf{degree function} {of $a(t)$} is defined by
	$$\delta_a(t) := \frac{2ta'(t)}{a(t)}.$$
	Furthermore, the \textbf{$k$-th degree function} {of $a(t)$} is defined by 
	$$\delta^{(k)}_a(t) := \frac{2t^k\left(\frac{\operatorname{d}^k}{\operatorname{d}t^k}a(t)\right)}{a(t)}.$$
	We say $a(t)$ has a \textbf{$k$-order finite degree}, if
	$$
\inf_{t\geqslant 0}\delta^{(i)}_a(t) = l^{(i)}_{a} > -\infty,~~~~\sup_{t\geqslant 0}\delta^{(i)}_a(t) = d^{(i)}_{a} < +\infty,
$$
for any $1\leqslant i \leqslant k.$
\end{definition}

Therefore, (\ref{eigenproblem}) can be formally rewritten as    
\begin{align}\label{defharm}
\Delta_\varphi u + u\Psi(u^2, |\nabla u|^2) = 0,
\end{align}
where $\varphi$-Laplacian operator $\Delta_\varphi u= \operatorname{div} (\varphi(|\nabla u|^2)\nabla u)$  and 
%$$\Psi(u^2, |\nabla u|^2) = \delta_{a}(u^2)\varphi(|\nabla u|^2)\frac{|\nabla u|^2}{u^2} + b(u^2)$$
\begin{align}\label{Psi}
	\Psi(s, t) = \frac{\delta_{a}(s)\varphi(t)t}{s} + b(s),~\quad b(s)= \frac{\lambda\psi(s)}{a(s)}.
\end{align}
The existence and regularity of (\ref{eigenproblem}) are discussed at the beginning of Section 2.

Inspired by Nash--Moser iteration technique used by X. Wang and L. Zhang \cite{wang2010local} and J. He, \textit{et al.} \cite{he2023gradient}, we establish the Cheng--Yau estimate of equation (\ref{defharm}) as follows.
\begin{theorem}\label{mainthm}
	Let $(M^n,g)$ be a complete Riemannian $n$-manifold with Ricci curvature bounded from below by $\operatorname{Ric} \geqslant -K$ {where} $K\geqslant0$, and let $u$ be a positive solution of (\ref{eigenproblem}) on the ball $B(o,2R)\subset M$. Suppose $\mathcal{F}(s,t) = a(s)\varphi(t)$, where $\varphi(t)$ and $a(s)$ have 2-order finite degrees, and  satisfy that
	\begin{equation}\tag{C1}
		\delta_{\varphi}(t) \geqslant l_{\varphi} >-1, \quad\text{ for any } t\in [0,\infty);
	\end{equation}
	\begin{equation}\tag{C2}
		\frac{\left(\delta_\varphi(t) + \delta_{a}(s)+1 \right)^2}{n-1} - 2\delta_{\varphi}'(t)t -2\delta_{a}'(s)s \geqslant \gamma >0, \quad\text{ for any } s,t\in [0,\infty) ;
	\end{equation}
	\begin{align}\label{psi2}\tag{C3}
		\varTheta := \sup_{\substack{s\geqslant 0,\\t\in \mathbb{R}^+ -  I}} \left(\frac{2\left(\delta_\varphi(t) + \delta_{a}(s)+1 \right)}{n-1} +\delta_{\varphi}(t) + \delta_{a}(s) - \delta_{\psi}(s) \right)^2  < \frac{4\gamma}{n-1},
	\end{align}
	where
	\begin{align}\label{psi1}
		I:= \left\{t>0: b(t)\left[\frac{2\left(\delta_\varphi(t) + \delta_{a}(s)+1 \right)}{n-1} + \delta_\varphi(t)+ \delta_{a}(s) - \delta_{\psi}(s)  \right] \geqslant 0, \text{ for each }s\geqslant 0\right\}.
	\end{align}
	Then, there exists a constant $C$ which depends only on $n$, $\gamma$, the bounds of $\delta^{(k)}_\varphi$ and $\delta^{(k)}_a$ ($k=1,2$), and $\delta_b$ , such that
	\begin{align}\label{esti}
		\frac{|\nabla u|}{u}\leqslant C\frac{1+\sqrt{K}R}{R}
	\end{align}
	on $B(o,R)$.
\end{theorem}

\begin{remark}
	When $a(s) \equiv 1$ thus $\delta_{a}(s) \equiv 0$, condition (C1)--(C3) is in the same form as  \cite{shen2025feasibilitynashmoseriterationchengyautype}. Hence, Theorem \ref{mainthm} generalizes the one in \cite{shen2025feasibilitynashmoseriterationchengyautype}.
\end{remark}

\begin{remark}
	If we remove the assumption that $\mathcal{F}(s,t) = a(t)\varphi(t)$, after a small revision in the the proof of Theorem \ref{mainthm}, we can formally derive the similar Cheng--Yau estimate for the general operator $Q$. However, it is important to remark that  such operator lacks the necessary regularity, even in $\{x\in M: |\nabla u|\geqslant\varepsilon\}$. We will discuss this case in a subsequent paper.
\end{remark}

Based on the aforementioned estimate, we show an interesting consequences for nonlinear eigenvalue problem.

\begin{theorem}\label{nepthm}
	Let $M$ be a complete and non-compact Riemannian manifold with non-negative Ricci curvature, if the nonlinear eigenvalue problem (\ref{eigenproblem}) satisfies condition (C1)--(C3) in Theorem \ref{mainthm}, then the non-vanishing eigenfunction $u$ must be unbounded with respect to eigenvalue $\lambda \neq 0$. 
\end{theorem}

This manuscript is arranged as follows. In Section 2, we derive a Bochner-type inequality  which is a necessary tool in {the} Moser's iteration. Then we prove the generalized Cheng--Yau gradient estimate in Section 3. In the last section, we discuss nonlinear eigenvalue problem and a specific example for $\Delta_{p}\left(\sum_{i=1}^{m}a_iu^{q_i}\right)= \lambda u^r$.

\section{Preliminary}

{We consider a (weak) positive  solution $u\in C^1(\Omega)\cap W^{1,\alpha}(\Omega)$ of equation (\ref{defharm}) over a bounded domain $\Omega\subset M$, that is, $a(u^2)\varphi(|\nabla u|^2)|\nabla u|, \psi(u^2)u\in L^1(\Omega)$ and}
\begin{align}
	{- \int_{\Omega}a(u^2)\varphi(|\nabla u|^2)\langle \nabla u, \nabla \phi \rangle + \int_{\Omega}\lambda\psi(u^2)u\phi = 0 }
\end{align}
for any $v\in C^\infty_0(\Omega)$. Note that $a(u^2)>0$ due to $u>0$ on $\Omega$. So we can replace $\phi$ by $\phi/a(u^2)$, then there exists a sequence $\{\phi_k\}\subset C_0^\infty(\Omega)$, such that,
$$\phi_k \to  \frac{\phi}{a(u^2)} \text{ a.e. in } \Omega$$
$$\nabla\phi_k \to  \frac{\nabla\phi}{a(u^2)} + \frac{2ua'(u^2)\nabla u}{a(u^2)} \text{ a.e. in } \Omega$$
Then
\begin{align}
	- \int_{\Omega}\varphi(|\nabla u|^2)\langle \nabla u, \nabla \phi \rangle + \int_{\Omega}\left(\frac{2u^2a'(u^2)}{a(u^2)}\cdot\frac{\varphi(|\nabla u|^2)|\nabla u|^2}{u^2} + \frac{\lambda \psi(u^2)}{a(u^2)}\right)u\phi = 0
\end{align}
for any $\phi\in C^\infty_0(\Omega)$. This shows that the positive solutions of (\ref{eigenproblem}) and (\ref{defharm}) are equivalent in the weak sense.  
{Let $M_\varepsilon := \{x\in M: |\nabla u|^2 (x) > \varepsilon/2\}$ for any fixed $\varepsilon >0$. Owing to 
$$\inf_{t\in [\frac{\varepsilon^2}{4}, C]}\varphi(t)>0,$$
for any fixed $C>0$, we have  $u\in W^{2,2}_{\text{loc}}(\Omega\cap M_\varepsilon)$ according to Remark 2.7 in \cite{cianchi2018second}. It should be noted that $\Delta_{\varphi}$ is uniformly elliptic only in the regular part $M_\varepsilon$, so that the weak solution $u$ is in fact smooth in $\Omega\cap M_\varepsilon$.

Denote $f := \log u$, $H := |\nabla u|^2$ and $\hat{H} := H/u^2 = |\nabla f|^2$. It is easy to check that
\begin{align}\label{u2f_nabla}
	\nabla H = {u^2}\left(\nabla \hat{H} + 2\hat{H}\nabla f \right)\quad
	\text{ and }\quad
	\Delta u = {u}\left(\Delta f + \hat{H} \right).
\end{align}
Then (\ref{defharm}) reduces to
\begin{align}\label{u2f_harmonic}
\begin{aligned}
\varphi(H)\Delta f = - {\varphi'(H)H}\frac{\left\langle\nabla \hat{H}, \nabla f\right\rangle}{\hat{H}} - \left({2\varphi'(H)H} + {\varphi(H)}\right)\hat{H}  - {\Psi(u^2,H)},
\end{aligned}
\end{align}
where $\Psi$ is defined in (\ref{Psi}).

We will derive the following Bochner-type inequality for the modified linearization operator  $\mathcal{L}_{\mathcal{W},\varphi}$ defined by
\begin{align}\label{linear}
\begin{aligned}
\mathcal{L}_{\mathcal{W},\varphi}(\eta) :=& \operatorname{div}\left({\varphi(H)}\mathcal{W}(\eta)\nabla\eta + {2\varphi'(H)} \mathcal{W}(\eta)\left\langle \nabla u, \nabla \eta \right\rangle \nabla u \right)/\mathcal{W}(\eta).
\end{aligned}
\end{align}

\begin{lemma}\label{btf_l}
     Any positive solution $u$ of  (\ref{defharm}) satisfies that 
	\begin{align}\label{weighted_linear2}
	\begin{aligned}
	\mathcal{L}_{\mathcal{W},\varphi}(\hat{H}) &\geqslant 2\varphi(H)\operatorname{Ric}(\nabla f) + \varphi(H)\left(\delta_{\varphi}(H)+1+{\delta_\mathcal{W}(\hat{H})}\right)\frac{|\nabla \hat{H}|^2}{2\hat{H}} \\
	&~~~~   + \left[\frac{2\left(\delta_\varphi(H) + \delta_{a}(u^2)+1 \right)^2}{n-1} - 4\delta_{\varphi}'(H)H -4\delta_{a}'(u^2)u^2\right]\varphi(H) \hat{H}^2\\
	&~~~~ + 2b(u^2)\left[\frac{2\left(\delta_\varphi(H) + \delta_{a}(u^2)+1 \right)}{n-1} + \delta_\varphi(H)-\delta_b(u^2) \right]\hat{H}\\
	&~~~~ + {\delta_\mathcal{W}(\hat{H})}\varphi'(H)H\left(\frac{\left\langle\nabla \hat{H}, \nabla f\right\rangle}{\hat{H}}\right)^2 + A(u^2,H)\varphi(H)\left\langle \nabla \hat{H}, \nabla f\right\rangle\\
	&~~~~ + \frac{2\varphi(H)}{n-1}\left((\delta_{\varphi}(H)+1)\frac{\left\langle\nabla \hat{H}, \nabla f\right\rangle}{2\hat{H}} + \frac{b(u^2)}{\varphi(H)} \right)^2,
	\end{aligned}
	\end{align}
	where
	\begin{align}\label{defA}
	\begin{aligned}
	A(u^2,H):=&\frac{2(\delta_{\varphi}(H)+1)\left(\delta_\varphi(H) + \delta_{a}(u^2)+1 \right)}{n-1}-2(\delta_{\varphi}(H)+1)\\ &-2\delta_{\varphi}'(H)H-\delta_a(u^2)\left(\delta_{\varphi}(H)+2\right).
	\end{aligned}	
	\end{align}
\end{lemma}

\begin{proof}
	By the definition of linearization in (\ref{linear}),
	\begin{align}\label{btf_e1}
	\begin{aligned}
	\mathcal{L}_{\mathcal{W},\varphi}(\hat{H}) 
	&= \varphi(H)\Delta \hat{H} + \varphi'(H)\left\langle\nabla H, \nabla \hat{H} \right\rangle + 2\varphi'(H)H\frac{\left\langle \nabla f, \nabla \hat{H} \right\rangle }{\hat{H}}\Delta f  \\
	&\quad+2\left\langle \nabla\left( \varphi'(H)H\frac{\left\langle \nabla f, \nabla \hat{H} \right\rangle }{\hat{H}} \right), \nabla f\right\rangle.\\
	&\quad+\delta_{\mathcal{W}}(\hat{H})\varphi(H)\frac{|\nabla\hat{H}|^2}{2\hat{H}} + \delta_{\mathcal{W}}(\hat{H})\varphi'(H)H\left(\frac{\left\langle\nabla f, \nabla \hat{H}\right\rangle}{\hat{H}}\right)^2.
	\end{aligned}
	\end{align}
	Utilizing (\ref{u2f_nabla}), (\ref{u2f_harmonic}) and the standard Bochner formula of Laplacian \cite{li2012geometric}, namely,
	\begin{align*}
	\frac{1}{2}\Delta \hat{H} = |\nabla^2 f|^2 + \left\langle \nabla\Delta f, \nabla f\right\rangle + \operatorname{Ric}(\nabla f),
	\end{align*}
	we can infer from (\ref{btf_e1}) that 
	\begin{align}\label{btf_e2}
	\begin{aligned}
	\mathcal{L}_{\mathcal{W},\varphi}(\hat{H}) 
	=& 2\varphi(H)\left( |\nabla^2 f|^2 + \left\langle \nabla\Delta f, \nabla f\right\rangle + \operatorname{Ric}(\nabla f)\right) + \varphi'(H)u^2|\nabla  \hat{H}|^2\\
	& + 2\varphi'(H) {H}\left\langle\nabla f, \nabla \hat{H} \right\rangle+ 2\varphi'(H)H\frac{\left\langle \nabla f, \nabla \hat{H} \right\rangle }{\hat{H}}\Delta f   \\
	&+\delta_{\mathcal{W}}(\hat{H})\varphi(H)\frac{|\nabla\hat{H}|^2}{2\hat{H}} + \delta_{\mathcal{W}}(\hat{H})\varphi'(H)H\left(\frac{\left\langle\nabla f, \nabla \hat{H}\right\rangle}{\hat{H}}\right)^2\\
	&- 2\left\langle \nabla\left(\left(\varphi(H) + 2\varphi'(H)H\right)\hat{H}\right), \nabla f\right\rangle- 2\left\langle \nabla\varphi(H)\Delta f, \nabla f\right\rangle\\
	&- 2\left\langle \nabla\Psi(u^2,H), \nabla f\right\rangle.
	\end{aligned}
	\end{align}
	Then, calculate directly the last two terms on the RHS of (\ref{btf_e2}) as follows.
	\begin{align}\label{btf_e21}
	\begin{aligned}
	- 2\langle \nabla\varphi(H)\Delta f, \nabla f\rangle 
	=& -2\varphi(H)\left\langle\nabla\Delta f , \nabla f\right\rangle- 2\varphi'(H)H\Delta f\frac{\left\langle \nabla \hat{H}, \nabla f  \right\rangle }{\hat{H}}\\
	&+ 4\varphi'(H)H\left( \frac{\delta_\varphi(H)}{2}{\left\langle\nabla \hat{H}, \nabla f\right\rangle} + \left(\delta_{\varphi}(H) + 1\right)\hat{H}^2\right)\\
	&+2\delta_{a}(u^2)\delta_{\varphi}(H)\varphi(H)\hat{H}^2 + 2b(u^2)\delta_{\varphi}(H)\hat{H},
	\end{aligned}
	\end{align}
	and
	\begin{align}\label{btf_e23}
	\begin{aligned}
	- 2\left\langle \nabla\Psi(u^2,H), \nabla f\right\rangle =& - 2\left\langle \nabla\left(\delta_{a}(u^2)\varphi(H)\hat{H}+b(u^2)\right), \nabla f\right\rangle \\
	=& -\left(4\delta'_{a}(u^2)u^2 + 2 \delta_{a}(u^2)\delta_{\varphi}(H) \right)\varphi(H)\hat{H}^2\\
	&-\delta_{a}(u^2)\left(2+\delta_{\varphi}(H)\right)\varphi(H)\left\langle \nabla \hat{H}, \nabla f  \right\rangle\\
	&-4b'(u^2)u^2\hat{H}.
	\end{aligned}
	\end{align}

	 For the remaining terms, we use the formulas in \cite{shen2025feasibilitynashmoseriterationchengyautype}:
	 \begin{align*}
	 	\begin{aligned}
	 	&- 2\left\langle \nabla\left(\left(\varphi(H) + 2\varphi'(H)H\right)\hat{H}\right), \nabla f\right\rangle\\
	 	=& -2\left(\varphi(H) + 2\varphi'(H)H\right)\langle \nabla\hat{H}, \nabla f\rangle -2\left(\varphi(H) + 2\varphi'(H)H\right)'H\langle \nabla \hat{H}, \nabla f\rangle\\
	 	& - 4\left(\varphi(H) + 2\varphi'(H)H\right)'H\hat{H}^2,
	 	\end{aligned}
	 \end{align*}
	 and
	\begin{align*}
	\begin{aligned}
	\left(\varphi(H) + 2\varphi'(H)H\right)'H 
	= \delta'_\varphi(H)H\varphi(H) + \frac{1}{2}\delta_\varphi(H)^2\varphi(H) + \frac{1}{2}\delta_\varphi(H)\varphi(H).
	\end{aligned}
	\end{align*}
	Hence (\ref{btf_e2}) becomes
	\begin{align}\label{btf_final}
	\begin{aligned}
	\mathcal{L}_{\mathcal{W},\varphi}(\hat{H}) &= 2\varphi(H)\left( |\nabla^2 f|^2  + \operatorname{Ric}(\nabla f)\right) + \varphi(H)(\delta_{\mathcal{W}}(\hat{H})+\delta_{\varphi}(H))\frac{|\nabla  \hat{H}|^2}{2\hat{H}} \\
	&\quad+ 2b(u^2)\left(\delta_\varphi(H)-\delta_b(u^2)\right)\hat{H}-\varphi(H)\left(4\delta_\varphi'(H)H + 4\delta'_{a}(u^2)u^2\right)\hat{H}^2\\
	&\quad - 2\varphi(H)\left(\delta_\varphi(H)+1+\delta_\varphi'(H)H + \delta_a(u^2)\left(\delta_{\varphi}(H)+2\right)\right)\left\langle \nabla \hat{H}, \nabla f\right\rangle.
	\end{aligned}
	\end{align}

	Next, we need to estimate the Hessian term $|\nabla^2f|^2$,  by choosing a local orthonormal frame $\{e_i\}$ with $e_1 = \nabla f/|\nabla f|$. Then
	\begin{align}\label{f11}
	f_{11} = \frac{\left\langle\nabla \hat{H}, \nabla f\right\rangle}{2\hat{H}}
	\quad\text{and}\quad
	\sum_{i=1}^{n}f_{1i}^2 = \frac{|\nabla \hat{H}|^2}{4\hat{H}}.
	\end{align}
	It is  immediately deduced from (\ref{u2f_harmonic}) and (\ref{f11}) that
	\begin{align}\label{u2f_harmonic_orth}
	\begin{aligned}
	\sum_{i=2}^{n}f_{ii} &= -f_{11} + \Delta f\\
	&=-\left(\delta_{\varphi}(H)+1\right)f_{11} - \left(\delta_{\varphi}(H)+1\right)\hat{H}- \frac{\Psi(u^2,H)}{\varphi(H)}.
	\end{aligned}
	\end{align}
	Therefore,
	\begin{align*}
	\begin{aligned}
	|\nabla^2 f|^2 &\geqslant \sum_{i=1}^{n}f_{1i}^2 + \sum_{i=2}^{n}f_{ii}^2 \\
	&\geqslant  \sum_{i=1}^{n}f_{1i}^2 + \frac{1}{n-1}\left(\sum_{i=2}^{n}f_{ii} \right)^2\\
	&\geqslant  \sum_{i=1}^{n}f_{1i}^2 + \frac{1}{n-1}\left(\left(\delta_{\varphi}(H)+1\right)f_{11} + \left(\delta_{\varphi}(H)+1\right)\hat{H}+ \frac{\Psi(u^2,H)}{\varphi(H)}\right)^2\\
	&\geqslant  \sum_{i=1}^{n}f_{1i}^2 + \frac{1}{n-1}\left(\left(\delta_{\varphi}(H)+1\right)f_{11} + \left(\delta_{\varphi}(H)+\delta_{a}(u^2)+1\right)\hat{H}+ \frac{b(u^2)}{\varphi(H)}\right)^2.
	\end{aligned}
	\end{align*}
	Due to (\ref{f11}) and (\ref{Psi}), it infers that
	\begin{align}\label{esti_hess}
	\begin{aligned}
	|\nabla^2 f|^2 \geqslant &\frac{|\nabla \hat{H}|^2}{4\hat{H}}  + \frac{\left(\delta_{\varphi}(H)+\delta_{a}(u^2)+1\right)^2}{n-1}\hat{H}^2 + \frac{1}{n-1}\left((\delta_{\varphi}(H)+1)\frac{\left\langle\nabla \hat{H}, \nabla f\right\rangle}{2\hat{H}} + \frac{b(u^2)}{\varphi(H)} \right)^2\\
	&+ \frac{\left(\delta_\varphi(H) +1 \right)\left(\delta_{\varphi}(H)+\delta_{a}(u^2)+1\right)}{n-1}\left\langle\nabla \hat{H}, \nabla f\right\rangle + \frac{2\left(\delta_{\varphi}(H)+\delta_{a}(u^2)+1\right)b(u^2)}{(n-1)\varphi(H)}\hat{H}.
	\end{aligned}
	\end{align}
	Substituting (\ref{esti_hess}) for $|\nabla^2 f|^2$ in (\ref{btf_final}) yields the  assertion.

\end{proof}

It is known that a Riemannian manifold with Ricci curvature bounded below has finite Sobolev constant. The following Sobolev inequality is necessary in the iteration process.
\begin{theorem}[\cite{SaloffCoste1992UniformlyEO}]\label{sblv_thm}
	For $n>2$,  let $(M^n,g)$ be a complete Riemannian $n$-manifold with Ricci curvature bounded from below by $\operatorname{Ric} \geqslant -K$ {for some} $K\geqslant0$, then there exists $C$, depending only on $n$, such that for ball $B(R) \subset M$ with radius $R$ and volume $V(R)$, we have for any $f \in C_0^{\infty}(B(R))$,
	$$
	\begin{gathered}\label{sblv_inq}
	\left(\int_{B(R)}|f|^{2 m} \right)^{1 / m} \leq e^{C(1+\sqrt{K}R)} V^{-2 / n} R^2\left(\int_{B(R)}\left(|\nabla f|^2+R^{-2}|f|^2\right) \right), \\
	\end{gathered}
	$$
	where $m=n /(n-2)$. Meanwhile, for $n \leq 2$, the above inequality holds with $n$ replaced by any fixed $n^{\prime}>2$.
\end{theorem}

\section{Cheng--Yau gradient estimates}
Let $(M,g)$ be a complete Riemannian manifold and $u$ be a positive local solution over a domain $\Omega$ containing $o\in M$ . Without loss of generality, we can assume $|\nabla u |(o)\neq 0$. Otherwise, Theorem \ref{mainthm} holds naturally.

Firstly, we choose $\mathcal{W} \equiv 1$ in Lemma \ref{btf_l}, so that $\delta_{\mathcal{W}} \equiv 0$. Dropping out the nonnegative terms on the RHS of (\ref{weighted_linear2}) and using conditions (C1), (C2) and   $\operatorname{Ric} \geqslant -K$, we deduce that
\begin{align}\label{btf}
\begin{aligned}
\mathcal{L}_{1,\varphi}(\hat{H}) 
&\geqslant -2K\varphi(H)\hat{H} + 2\gamma\varphi(H) \hat{H}^2 - a_0\varphi(H)|\nabla \hat{H}||\nabla f|\\
&~~~~ +
2b(u^2)\left[\frac{2\left(\delta_\varphi(H) + \delta_a(u^2) +1 \right)}{n-1} + \delta_\varphi(H)-\delta_b(u^2) \right]\hat{H},
\end{aligned}
\end{align}
where $|A(u^2, H)| \leqslant a_0$, since the 2-order finite degree functions of $a(s)$ and $\varphi(t)$. We will henceforth use $a_i$ to denote all  the constants related to $n$, $\varphi$, $a$ and $b$.

Since the condition (C3) and the identity $\delta_b = \delta_\psi-\delta_a$ infer $$b(u^2)\left[\frac{2\left(\delta_\varphi(H) + \delta_{a}(u^2)+1 \right)}{n-1} + \delta_\varphi(H)-\delta_b(u^2) \right]\geqslant0$$
on $\{u^2\in I\}$, we take the  
 test function $\phi$ compactly supported in ${\Omega\cap M_\varepsilon}$, defiend by
$$\phi :=\frac{\chi\hat{H}_\varepsilon^{q}\eta^2}{\varphi(H)}, $$ 
 where $\chi(x)$ is the characteristic function of $\left\{x\in \Omega: u^2(x)\in I\right\}$, $\hat{H}_\varepsilon := \left(\hat{H} -\varepsilon\right)^+$. The cutoff function $\eta\in C^\infty_0(\Omega)$ and constant $q>1$ will be determined later.
Then, it follows that
\begin{align}\label{wbtf}
\begin{aligned}
&\int_{{\Omega\cap M_\varepsilon}} \left\langle \varphi(H)\nabla \hat{H} + 2\varphi'(H) \left\langle \nabla u, \nabla \hat{H}   \right\rangle \nabla u, \nabla \phi \right\rangle
\\&\leqslant 2K\int_{{\Omega\cap M_\varepsilon}}\varphi(H)\hat{H}\phi - 2\gamma\int_{{\Omega\cap M_\varepsilon}} \varphi(H) \hat{H}^2\phi +  a_0\int_{{\Omega\cap M_\varepsilon}}\varphi(H)|\nabla \hat{H}||\nabla f|\phi.
\end{aligned}
\end{align}

Substituting the gradient of this test function
\begin{align*}
\begin{aligned}
\nabla\phi &= \frac{q\chi\hat{H}_\varepsilon^{q-1}\eta^2}{\varphi(H)} \nabla\hat{H} + \frac{2\chi\hat{H}_\varepsilon^q\eta}{\varphi(H)}\nabla \eta - \frac{\chi\varphi'(H)\hat{H}_\varepsilon^{q}\eta^2}{\varphi(H)^2}\nabla H \\
&= \left(\frac{q}{\varphi(H)}\hat{H}_\varepsilon^{q-1}\chi - \frac{\varphi'(H)H}{\varphi(H)^2}\frac{\hat{H}_\varepsilon^{q}\chi}{\hat{H}} \right)\eta^2\nabla\hat{H} + \frac{2\hat{H}_\varepsilon^q\eta\chi}{\varphi(H)}\nabla \eta -\frac{2\varphi'(H)H}{\varphi(H)^2}{\hat{H}_\varepsilon^{q}\chi}\nabla f.
\end{aligned}
\end{align*}
to the LHS of (\ref{wbtf}), we have
\begin{align}\label{wbtf_lhs}
\begin{aligned}
\int_{{\Omega\cap M_\varepsilon}}  &\left\langle \varphi(H)\nabla \hat{H} + 2\varphi'(H) \left\langle \nabla u, \nabla \hat{H}   \right\rangle \nabla u, \nabla \phi \right\rangle
\\
=&\int_{{\Omega\cap M_\varepsilon}}\left({q}\hat{H}_\varepsilon^{q-1} - \frac{\delta_{\varphi}(H)}{2}\frac{\hat{H}_\varepsilon^{q}}{\hat{H}} \right)\left(|\nabla \hat{H}|^2 + \frac{2\varphi'(H)}{\varphi(H)} \left\langle \nabla u, \nabla \hat{H}\right\rangle^2 \right)\chi\eta^2\\
&+2\int_{{\Omega\cap M_\varepsilon}}\hat{H}_\varepsilon^q \left(\left\langle \nabla \hat{H}, \nabla \eta\right\rangle + \frac{2\varphi'(H)}{\varphi(H)}\left\langle \nabla u, \nabla \hat{H}\right\rangle\left\langle \nabla u, \nabla \eta\right\rangle\right)\chi\eta\\
&-\int_{{\Omega\cap M_\varepsilon}}\delta_{\varphi}(H)\hat{H}_\varepsilon^{q}\left(\left\langle \nabla f, \nabla \hat{H}\right\rangle + \frac{2\varphi'(H)}{\varphi(H)} \left\langle \nabla u, \nabla \hat{H}\right\rangle\left\langle \nabla u, \nabla f\right\rangle\right)\chi\eta^2.
\end{aligned}
\end{align}
Because $\hat{H}_\epsilon \leqslant \hat{H}$, the first term on the RHS of (\ref{wbtf_lhs}) is estimated by
\begin{align}\label{wbtf_lhs1}
\begin{aligned}
\int_{{\Omega\cap M_\varepsilon}}&\left({q}\hat{H}_\varepsilon^{q-1} - \frac{\delta_{\varphi}(H)}{2}\frac{\hat{H}_\varepsilon^{q}}{\hat{H}} \right)\left(|\nabla \hat{H}|^2 + \frac{2\varphi'(H)}{\varphi(H)} \left\langle \nabla u, \nabla \hat{H}\right\rangle^2 \right)\chi\eta^2\\
\geqslant& \int_{\{\varphi'(H) < 0\}}q\hat{H}_\varepsilon^{q-1}\left(|\nabla \hat{H}|^2 + l_\varphi|\nabla \hat{H}|^2 \right)\chi\eta^2 + \int_{\{\varphi'(H) \geqslant 0\}}\left({q} - \frac{d_\varphi}{2}\right)\hat{H}_\varepsilon^{q-1}|\nabla \hat{H}|^2 \chi\eta^2\\
\geqslant& {a_1q} \int_{{\Omega\cap M_\varepsilon}}\hat{H}_\varepsilon^{q-1}|\nabla \hat{H}|^2\chi\eta^2,
\end{aligned}
\end{align}
if we choose $q>d_\varphi$, the upper bound of $\delta_{\varphi}$. Moreover, the second term on the RHS of (\ref{wbtf_lhs}) is 
\begin{align}\label{wbtf_lhs2}
\begin{aligned}
&2\int_{{\Omega\cap M_\varepsilon}}\hat{H}_\varepsilon^q \left(\left\langle \nabla \hat{H}, \nabla \eta\right\rangle + \frac{2\varphi'(H)}{\varphi(H)}\left\langle \nabla u, \nabla \hat{H}\right\rangle\left\langle \nabla u, \nabla \eta\right\rangle\right)\chi\eta \\
\geqslant& -2\int_{{\Omega\cap M_\varepsilon}}\hat{H}_\varepsilon^q \left(|\nabla \hat{H}|| \nabla \eta| + \left|\frac{2\varphi'(H)}{\varphi(H)}\right||\nabla u|^2 |\nabla \hat{H}||\nabla \eta|\right)\chi\eta\\
\geqslant&-a_2 \int_{{\Omega\cap M_\varepsilon}}\hat{H}_\varepsilon^q |\nabla \hat{H}|| \nabla \eta|\chi\eta.
\end{aligned}
\end{align}
Finally, we estimate the last term on the RHS of (\ref{wbtf_lhs})
\begin{align}\label{wbtf_lhs3}
\begin{aligned}
&-\int_{{\Omega\cap M_\varepsilon}}\delta_{\varphi}(H)\hat{H}_\varepsilon^{q}\left(\left\langle \nabla f, \nabla \hat{H}\right\rangle + \frac{2\varphi'(H)}{\varphi(H)} \left\langle \nabla u, \nabla \hat{H}\right\rangle\left\langle \nabla u, \nabla f\right\rangle\right)\chi\eta^2\\
\geqslant& -\int_{{\Omega\cap M_\varepsilon}}\left|\delta_{\varphi}(H)\right|\hat{H}_\varepsilon^{q}\left(|\nabla f||\nabla \hat{H}| + \left|\delta_{\varphi}(H)\right| |\nabla f||\nabla \hat{H}|\right)\chi\eta^2\\
\geqslant& -a_3\int_{{\Omega\cap M_\varepsilon}}\hat{H}_\varepsilon^{q}|\nabla f||\nabla \hat{H}|\chi\eta^2.
\end{aligned}
\end{align}

Combining (\ref{wbtf_lhs1}), (\ref{wbtf_lhs2}) (\ref{wbtf_lhs3}) with (\ref{wbtf_lhs}), then (\ref{wbtf}) leads to
\begin{align*}
\begin{aligned}
{a_1 q}\int_{{\Omega\cap M_\varepsilon}}&\hat{H}_\varepsilon^{q-1}|\nabla \hat{H}|^2\chi\eta^2 + 2\gamma\int_{{\Omega\cap M_\varepsilon}} \hat{H}^2\hat{H}_\varepsilon^{q}\chi\eta^2\\ \leqslant& 2K\int_{{\Omega\cap M_\varepsilon}}\hat{H}\hat{H}_\varepsilon^{q}\chi\eta^2  +  (a_0 + a_3)\int_{{\Omega\cap M_\varepsilon}}\hat{H}_\varepsilon^{q}|\nabla \hat{H}||\nabla f|\chi\eta^2
+ a_2\int_{{\Omega\cap M_\varepsilon}}\hat{H}_\varepsilon^{q}|\nabla \hat{H}||\nabla \eta|\chi\eta
\\ \leqslant& 2K\int_{\Omega}\hat{H}^{q+1}\chi\eta^2  +  (a_0 + a_3)\int_{\Omega}\hat{H}^{q}|\nabla \hat{H}||\nabla f|\chi\eta^2
+ a_2\int_{\Omega}\hat{H}^{q}|\nabla \hat{H}||\nabla \eta|\chi\eta.
\end{aligned}.
\end{align*}
 Due to the compact support of $\hat{H}_\varepsilon\eta$ in $ M_\varepsilon\cap\Omega$, after extending the integrals to $\Omega$,  we obtain by Fatou's lemma that
\begin{align}\label{int_esti_1}
\begin{aligned}
{a_1 q} \int_{\Omega}&\hat{H}^{q-1}|\nabla \hat{H}|^2\chi\eta^2 + 2\gamma\int_{\Omega} \hat{H}^{q+2}\chi\eta^2\\ \leqslant&{\lim\limits_{\varepsilon \to 0}{a_1 q} \int_{\Omega}\hat{H}_\varepsilon^{q-1}|\nabla \hat{H}|^2\chi\eta^2 + \lim\limits_{\varepsilon \to 0}2\gamma\int_{\Omega} \hat{H}^2\hat{H}_\varepsilon^{q}\chi\eta^2}
\\ \leqslant& 2K\int_{\Omega}\hat{H}^{q+1}\chi\eta^2  +  (a_0 + a_3)\int_{\Omega}\hat{H}^{q}|\nabla \hat{H}||\nabla f|\chi\eta^2
+ a_2\int_{\Omega}\hat{H}^{q}|\nabla \hat{H}||\nabla \eta|\chi\eta.
\end{aligned}
\end{align}

By {Cauchy's} inequality, the last two terms on the RHS of (\ref{int_esti_1}) could be estimated by
\begin{align*}
\begin{aligned}
(a_0 + a_3)\int_{\Omega}\hat{H}^{q}|\nabla \hat{H}||\nabla f|\chi\eta^2 \leqslant \frac{(a_0 + a_3)^2}{4\gamma}\int_{\Omega} \hat{H}^{q-1}|\nabla \hat{H}|^2\chi\eta^2+ \gamma\int_{\Omega} \hat{H}^{q+2}\chi\eta^2,
\end{aligned}
\end{align*}
and
\begin{align*}
\begin{aligned}
a_2\int_{\Omega}\hat{H}^{q}|\nabla \hat{H}||\nabla \eta|\chi\eta \leqslant \frac{a_1q}{2}\int_{\Omega} \hat{H}^{q-1}|\nabla \hat{H}|^2\chi\eta^2+ \frac{a_2^2}{2a_1q}\int_{\Omega} \hat{H}^{q+1}|\nabla \eta|^2\chi.
\end{aligned}
\end{align*}
Then, (\ref{int_esti_1}) becomes
\begin{align}\label{int_esti_2}
\begin{aligned}
\frac{a_1 q}{4} \int_{\Omega}  \hat{H}^{q-1}|\nabla \hat{H}|^2\chi\eta^2 + \gamma\int_{\Omega} \hat{H}^{q+2}\chi\eta^2
\leqslant 2K\int_{\Omega} \hat{H}^{q+1}\chi\eta^2
+ \frac{a_2^2}{2a_1q}\int_{\Omega} \hat{H}^{q+1}|\nabla \eta|^2\chi.
\end{aligned}
\end{align}
where we choose $q$ large enough such that
\begin{align}\label{require_b}
q>\max\left\{\frac{(a_0 + a_3)^2}{a_1\gamma},d_{\varphi},1\right\}.
\end{align}

Since
\begin{align*}
\left|\nabla\left(\hat{H}^{q/2+1/2}\eta\right)\right|^2 &\leqslant \frac{1}{2}\left(q+1\right)^2\hat{H}^{q-1}|\nabla \hat{H}|^2\eta^2 + 2\hat{H}^{q+1}|\nabla \eta|^2 \\
&\leqslant 2q^2\hat{H}^{q-1}|\nabla \hat{H}|^2\eta^2 + 2\hat{H}^{q+1}|\nabla \eta|^2,
\end{align*}
it follows that 
\begin{align}\label{int_esti}
\begin{aligned}
\int_{\Omega} \left|\nabla\left(\hat{H}^{q/2+1/2}\eta\right)\right|^2\chi + \frac{8\gamma q}{a_1} \int_{\Omega} \hat{H}^{ q+2}\chi\eta^2
\leqslant \frac{16Kq}{a_1}\int_{\Omega} \hat{H}^{q+1}\chi\eta^2
+ a_3\int_{\Omega} \hat{H}^{q+1}|\nabla \eta|^2\chi.
\end{aligned}
\end{align}

On the other hand, to deal with the part that $u^2 \notin I$, we set $\mathcal{W}(\eta) = \eta^\alpha$ for some $\alpha>0$ which will be given later, then $\delta_\mathcal{W} \equiv 2\alpha$ and (\ref{weighted_linear2}) becomes
\begin{align}\label{wt}
\begin{aligned}
\frac{\mathcal{L}_{\alpha,\varphi}(\hat{H})}{\varphi(H)} &\geqslant -2K\hat{H}   + 2\gamma\hat{H}^2 + \left(\delta_{\varphi}(H)+1+{2\alpha}\right)\frac{|\nabla \hat{H}|^2}{2\hat{H}}\\
&\quad - a_0|\nabla \hat{H}||\nabla f|+\left[\frac{(\delta_{\varphi}(H)+1)^2}{2(n-1)} + {\alpha}\delta_{\varphi}(H)\right]\frac{\left\langle\nabla \hat{H}, \nabla f\right\rangle^2}{\hat{H}^2}\\
&\quad  + \frac{2(\delta_{\varphi}(H)+1)}{(n-1)}\cdot\frac{\left\langle\nabla \hat{H}, \nabla f\right\rangle}{\hat{H}}\frac{b(u^2)}{\varphi(H)}+ \frac{2}{n-1}\cdot\frac{b(u^2)^2}{\varphi(H)^2}\\
&\quad  + \frac{2b(u^2)}{\varphi(H)}\left[\frac{2\left(\delta_\varphi(H)+\delta_{a}(u^2) +1 \right)}{n-1} + \delta_\varphi(H)-\delta_b(u^2) \right]\hat{H}\\
&\geqslant -2K\hat{H}   + 2\gamma\hat{H}^2 - a_0|\nabla \hat{H}||\nabla f|\\
&\quad +\left[\frac{(\delta_{\varphi}(H)+1)^2}{2(n-1)} + {\alpha}\left(\delta_{\varphi}(H)+1 \right)\right]\frac{\left\langle\nabla \hat{H}, \nabla f\right\rangle^2}{\hat{H}^2}\\
&\quad  + \frac{2(\delta_{\varphi}(H)+1)}{(n-1)}\cdot\frac{\left\langle\nabla \hat{H}, \nabla f\right\rangle}{\hat{H}}\frac{b(u^2)}{\varphi(H)}+ \frac{2}{n-1}\cdot\frac{b(u^2)^2}{\varphi(H)^2} + 2B(u^2,H)\hat{H}\frac{b(u^2)}{\varphi(H)}, 
\end{aligned}
\end{align}
where the last inequality follows from
$$
\left(\delta_{\varphi}(H)+1+{2\alpha}\right)\frac{|\nabla \hat{H}|^2}{2\hat{H}}\geqslant\frac{\alpha|\nabla \hat{H}|^2}{\hat{H}}\geqslant\frac{\alpha\left\langle\nabla \hat{H}, \nabla f\right\rangle^2}{\hat{H}^2},
$$
and we denote
$$B(s,t):= \frac{2\left(\delta_\varphi(t) + \delta_{a}(s) + 1 \right)}{n-1} + \delta_\varphi(t)-\delta_b(s)$$
for short. By using the inequality $x^2 + 2xy \geqslant -y^2$ twice, we have
\begin{align*}
\begin{aligned}
\frac{\mathcal{L}_{\alpha,\varphi}(\hat{H})}{\varphi(H)}
&\geqslant -2K\hat{H}   + 2\gamma\hat{H}^2- a_0|\nabla \hat{H}||\nabla f|\\
&~~~~ +\left[\frac{2}{n-1} - \frac{2(\delta_{\varphi}(H)+1)^2}{2\alpha(\delta_{\varphi}(H)+1)(n-1)^2+(\delta_{\varphi}(H)+1)^2(n-1) }\right]\frac{b(u^2)^2}{\varphi(H)^2} \\
&~~~~ +  2B(u^2,H)\hat{H}\cdot\frac{b(u^2)}{\varphi(H)}\\
&\geqslant -2K\hat{H} - a_0|\nabla \hat{H}||\nabla f| +2\left[\gamma-B(u^2,H)^2\left(\frac{n-1}{4}+\frac{(\delta_{\varphi}(H)+1)}{8\alpha}\right)\right]\hat{H}^2\\
&\geqslant -2K\hat{H} - a_0|\nabla \hat{H}||\nabla f| +2\left[\gamma-B(u^2,H)^2\left(\frac{n-1}{4}+\frac{(d_{\varphi}+1)}{8\alpha}\right)\right]\hat{H}^2.
\end{aligned}
\end{align*}

Similarly, choose the test function
$$\phi :=\frac{\bar{\chi}\hat{H}_\varepsilon^{q-\alpha}\eta^2}{\varphi(H)}, $$
where $\bar{\chi}$ is the characteristic function of $\left\{x\in \Omega: u^2(x)\notin I\right\}$, then,
\begin{align}\label{wbtf2}
\begin{aligned}
&\int_{M_\varepsilon} \left\langle \varphi(H)\hat{H}^\alpha\nabla \hat{H} + 2\varphi'(H)\hat{H}^\alpha \left\langle \nabla u, \nabla \hat{H}   \right\rangle \nabla u, \nabla \phi \right\rangle
\\\leqslant& 2K\int_{M_\varepsilon}\varphi(H)\hat{H}^{1+\alpha}\phi  +  a_0\int_{M_\varepsilon}\varphi(H)|\nabla \hat{H}||\nabla f|\hat{H}^{\alpha}\phi\\
&-
2\int_{M_\varepsilon}\left[\gamma-B(u^2,H)^2\left(\frac{n-1}{4}+\frac{(d_{\varphi}+1)}{8\alpha}\right)\right]\varphi(H)\hat{H}^{2+\alpha}\phi.
\end{aligned}
\end{align}

Now we set
$$\alpha(\gamma,d_{\varphi},\varTheta) :=\frac{(d_\varphi+1)}{4\theta}\sup_{\substack{s>0,\\t\in \mathbb{R}^+ -  I}} B(s,t)^2,$$ 
where the constant
$$\theta(\gamma,\varTheta):= \gamma-\frac{n-1}{4}\sup_{\substack{s>0,\\t\in \mathbb{R}^+ -  I}} B(s,t)^2 $$
must be positive due to the condition (\ref{psi2}). Then the last term on the RHS of (\ref{wbtf2}) could be estimated from below as
\begin{align*}
\begin{aligned}
&2\int_{M_\varepsilon}\left[\gamma-B(u^2,H)^2\left(\frac{n-1}{4}+\frac{(d_{\varphi}+1)}{8\alpha}\right)\right]\varphi(H)\hat{H}^{2+\alpha}\phi\\
=&2\int_{\{u^2(x)\notin I\}}\left[\gamma-B(u^2,H)^2\left(\frac{n-1}{4}+\frac{(d_{\varphi}+1)}{8\alpha}\right)\right]\varphi(H)\hat{H}^{2+\alpha}\phi\\
\geqslant&2\int_{\{u^2(x)\notin I\}}\left[\gamma-\sup_{\substack{s>0,\\t\in \mathbb{R}^+ -  I}}B(s,t)^2\left(\frac{n-1}{4}+\frac{(d_{\varphi}+1)}{8\alpha}\right)\right]\varphi(H)\hat{H}^{2+\alpha}\phi\\
\geqslant&2\int_{\{u^2(x)\notin I\}}\left[\theta-\sup_{\substack{s>0,\\t\in \mathbb{R}^+ -  I}}B(s,t)^2\left(\frac{d_{\varphi}+1}{8\alpha}\right)\right]\varphi(H)\hat{H}^{2+\alpha}\phi\\
=&\theta\int_{M_\varepsilon}\hat{H}^{2+\alpha}\hat{H}_\varepsilon^{q-\alpha}\eta^2\bar{\chi}.
\end{aligned}
\end{align*}
Following the same process from (\ref{wbtf_lhs}) to (\ref{int_esti}), we can also obtain
\begin{align}\label{int_esti_c}
\begin{aligned}
\int_{\Omega} \left|\nabla\left(\hat{H}^{q/2+1/2}\eta\right)\right|^2\bar{\chi} + &a_4\theta q \int_{\Omega} \hat{H}^{ q+2}\bar{\chi}\eta^2\\
\leqslant& a_5Kq\int_{\Omega} \hat{H}^{q+1}\bar{\chi}\eta^2
+ a_6\int_{\Omega} \hat{H}^{q+1}|\nabla \eta|^2\bar{\chi},
\end{aligned}
\end{align}
for constant 
\begin{align}\label{require_b2}
q>\max\left\{\frac{a_7}{\theta},d_{\varphi},\alpha\right\}.
\end{align}
Noticing that $\chi + \bar{\chi} \equiv 1$, we can adjust the constants and combine (\ref{int_esti}) and (\ref{int_esti_c}) to get
\begin{align}\label{int_esti_f}
\begin{aligned}
\int_{\Omega} \left|\nabla\left(\hat{H}^{q/2+1/2}\eta\right)\right|^2 + &a_7 q \int_{\Omega} \hat{H}^{ q+2}\eta^2\\
\leqslant& a_8Kb\int_{\Omega} \hat{H}^{q+1}\eta^2
+ a_9\int_{\Omega} \hat{H}^{q+1}|\nabla \eta|^2.
\end{aligned}
\end{align}

Next we set the domain $\Omega$ to be the geodesic ball $B(o,R)$ (also denoted by $B(R)$ for short). Theorem \ref{sblv_thm} shows that 
\begin{align}
\begin{aligned}\label{int_esti_sblv}
\left(\int_{\Omega} \hat{H}^{(q+1)m} \eta^{2m}\right)^{1/m}\leqslant e^{C(1+\sqrt{K}R)}V^{-2/n}\left(R^2\int_{\Omega}\left|\nabla\left(\hat{H}^{q/2+1/2}\eta\right)\right|^2+ \int_{\Omega} \hat{H}^{q+1} \eta^{2}\right).
\end{aligned}
\end{align}
with $m=n/(n-2)$, when $n>2$. Set $q_0 = c_0\left(1+\sqrt{K}R\right)$ and choose $c_0$ large enough to satisfy (\ref{require_b}) and (\ref{require_b2}). Apply (\ref{int_esti_sblv}) to (\ref{int_esti_c}), then direct calculation implies that
\begin{align*}
\begin{aligned}
&\left(\int_{\Omega} \hat{H}^{(q+1)m} \eta^{2m}\right)^{1/m} + a_7qe^{c_1q_0}\left(\frac{R^2}{V^{2/n}}\right)\int_{\Omega} \hat{H}^{q+2}\eta^2\\
&~~~~~~~~~~~~~\leqslant (a_{8}KR^2q+1)e^{c_1q_0}V^{-2/n}\int_{\Omega} \hat{H}^{q+1}\eta^{2} + a_{9}e^{c_1q_0}\left(\frac{R^2}{V^{2/n}}\right)\int_{\Omega} \hat{H}^{q+1}|\nabla\eta|^2.
\end{aligned}
\end{align*}
Noticing that
$$a_{8}KR^2q+1 \leqslant \max\left\{{a_8},1\right\}\cdot(KR^2 + 1)q \leqslant \max\left\{{a_8},1\right\}\cdot(\sqrt{K}R + 1)^2q,$$
we have
\begin{align}\label{int_esti_final}
\begin{aligned}
&\left(\int_{\Omega} \hat{H}^{(q+1)m} \eta^{2m}\right)^{1/m} + a_7qe^{c_1q_0}\left(\frac{R^2}{V^{2/n}}\right)\int_{\Omega} \hat{H}^{q+2}\eta^2\\
&~~~~~~~~~~~~~\leqslant a_{10}q_0^2qe^{c_1q_0}V^{-2/n}\int_{\Omega} \hat{H}^{q+1}\eta^{2} + a_{10}e^{c_1q_0}\left(\frac{R^2}{V^{2/n}}\right)\int_{\Omega} \hat{H}^{q+1}|\nabla\eta|^2.
\end{aligned}
\end{align}

In the following lemma, we will make the most use of the high-order term ($\hat{H}^{q+2}$) to show an $L^{\beta}$ bound of $\hat{H}$, so that we can let such $\beta$ to be the initial value of the iteration.
\begin{lemma}\label{Lb_esti}
	Under the same conditions above, take
	\begin{align}
	\beta_0 = q_0+1\text{ and } \beta_1 = \beta_0m.
	\end{align}
	Then there exists $a_{11}>0$ such that
	\begin{align}
	||\hat{H}||_{L^{\beta_1}(B{(3R/4)})} \leqslant a_{11}\frac{q^2_0}{R^2}V^{1/\beta_1}.
	\end{align}
\end{lemma}
\begin{proof}
	Set $q = q_0$ and decompose the first term on RHS of (\ref{int_esti_final}) into two parts as
	\begin{align*}
	\begin{aligned}
	a_{10}q_0^3e^{c_1q_0}V^{-2/n}\int_{B(R)} \hat{H}^{\beta_0}\eta^{2} = a_{10}q_0^3e^{c_1q_0}V^{-2/n}\left(\int_{\Omega_1} \hat{H}^{\beta_0}\eta^{2} + \int_{\Omega_2} \hat{H}^{\beta_0}\eta^{2}\right),
	\end{aligned}
	\end{align*}
	where
	\begin{align*}
	\Omega_1:= \left\{\hat{H}>\frac{2a_{10}q_0^2}{a_7R^2}\right\} \text{ and } \Omega_2:= \left\{\hat{H}\leqslant\frac{2a_{10}q_0^2}{a_7R^2}\right\}.
	\end{align*}
	This yields
	\begin{align}\label{RHS1_esti}
	\begin{aligned}
	a_{10}q_0^3e^{c_1q_0}V^{-2/n}&\int_{B_R} \hat{H}^{\beta_0}\eta^{2} \\
	\leqslant& \frac{1}{2}a_7q_0e^{c_1q_0}\left(\frac{R^2}{V^{2/n}}\right)\int_{B_R} \hat{H}^{\beta_0+1}\eta^2 +
	a_{10}a_{12}^{\beta_0}q_0^3e^{c_1q_0}V^{1-2/n}\left(\frac{q_0}{R}\right)^{2\beta_0}.
	\end{aligned}
	\end{align}

	To deal with the second term on the RHS of (\ref{int_esti_final}), we should choose the cutoff function $\eta = \eta_0^{\beta_0+1}$, where $\eta_0$ is a smooth function with compact support in $B(R)$ satisfying $0\leqslant \eta_0 \leqslant 1$, $\eta_0 \equiv 1$ in $B{(3R/4)}$ and 
	$
	|\nabla \eta_0|\leqslant{c_2(n)}/{R}.
	$
	Therefore,
	\begin{align}\label{eta_esti}
	|\nabla \eta|^2\leqslant c^2_2\left(\frac{\beta_0+1}{R} \right)^2\eta_0^{2\beta_0} = c^2_2\left(\frac{\beta_0+1}{R}\right)^2\eta^{\frac{2\beta_0}{\beta_0+1}}.
	\end{align}
	Substituting (\ref{eta_esti}) for the second term on the RHS of (\ref{int_esti_final}), we have
	\begin{align}\label{RHS2_esti}
	\begin{aligned}
	a_{10}q_0&  e^{c_1q_0}\left(\frac{R^2}{V^{2/n}}\right)\int_{B(R)} \hat{H}^{\beta_0}|\nabla\eta|^2 \leqslant a_{10}c_2^2\left(\beta_0+1\right)^2e^{c_1q_0}V^{-2/n}\int_{B(R)} \hat{H}^{\beta_0}\eta^{\frac{2\beta_0}{\beta_0+1}}\\
	\leqslant& a_{10}c_2^2e^{c_1q_0}\left(\frac{\left(\beta_0+1\right)^2}{V^{2/n}}\right)\left(\int_{B(R)} \hat{H}^{\beta_0+1}\eta^{2}\right)^{\frac{\beta_0}{\beta_0+1}}\left(\int_{B(R)} 1\right)^{\frac{1}{\beta_0+1}}\\
	\leqslant&\frac{1}{2}a_7q_0e^{c_1q_0}\left(\frac{R^2}{V^{2/n}}\right)\int_{B(R)} \hat{H}^{\beta_0+1}\eta^{2}+ 2a_{10}c_2^2\left(\frac{4a_{10}c_2^2}{a_7}\right)^{\beta_0}\beta_0e^{c_1q_0}V^{1-2/n}\left(\frac{\beta_0}{R}\right)^{2\beta_0},
	\end{aligned}
	\end{align}
	where we use the H\"older's inequality and Young's inequality at the last two inequalities, respectively.

	It follows from (\ref{int_esti_final}), (\ref{RHS1_esti}), and (\ref{RHS2_esti})  that
	\begin{align}
	\begin{aligned}
	\left(\int_{B(R)} \hat{H}^{\beta_0m} \eta^{2m}\right)^{1/m} \leqslant a_7a_{8}^{\beta_0}\beta_0^3e^{c_1q_0}V^{1/m}\left(\frac{\beta_0}{R}\right)^{2\beta_0},
	\end{aligned}
	\end{align}
	which is exact (\ref{Lb_esti}) after taking $(1/\beta_0)$-root on both sides.
\end{proof}

Now, we start to prove our main theorem.

\begin{proof}[Proof of Theorem \ref{mainthm}]
	Here we go back to (\ref{int_esti_final}) and dismiss the second nonnegative term on the LHS. It follows that
	\begin{align}\label{int_esti_iteration}
	\begin{aligned}
	\left(\int_{B(R)} \hat{H}^{(q+1)m} \eta^{2m}\right)^{1/m}
	\leqslant a_{10}\left(\frac{e^{c_1q_0}}{V^{2/n}}\right)\int_{B(R)} \left(\beta_0^2(q+1)\eta^{2}+ R^2|\nabla\eta|^2 \right)\hat{H}^{q+1},
	\end{aligned}
	\end{align}
	where $\beta_0$ are given in Lemma \ref{Lb_esti}.
	
	We now choose the sequences $\beta_k = \beta_0m^k$, $q_k+1 = \beta_k$, and  radii $R_k = \frac{R}{2} + \frac{R}{4^k}$, which yield a sequence of
	cutoff functions $\eta_k$ such that
	\begin{align}
	\left\{\begin{aligned}
	\eta_k \equiv 1~~~~~~~~~~~~~~~~~~ &\text{ ~~~in } B(R_{k+1}),\\
	0 \leqslant \eta_k \leqslant 1 \text{ and } |\nabla\eta_k|\leqslant\frac{c_3(n)4^k}{R} &\text{ ~~~in } B(R_{k}) - B(R_{k+1}),\\
	\eta_k \equiv 0~~~~~~~~~~~~~~~~~~ &\text{ ~~~in } B(R) - B(R_{k+1}).
	\end{aligned}\right.
	\end{align}
	Letting $q = q_k$ in (\ref{int_esti_iteration})and noting that
	$
	q_k < \beta_{k} = q_k + 1 \leqslant 2q_k,
	$
	we have
	\begin{align*}
	\begin{aligned}
	\left(\int_{B(R_{k+1})} \hat{H}^{\beta_{k+1}}\right)^{1/m}
	&\leqslant \left(\frac{a_{10}e^{c_1q_0}}{V^{2/n}}\right)\left(\beta_0^3m^k+ c_3 16^k \right)\int_{B(R_k)} \hat{H}^{\beta_k}\\
	&\leqslant \left(\frac{a_{10}\left(\beta_0^3+c_3\right)e^{c_1q_0}}{V^{2/n}}\right)\left( 16 \right)^k\int_{B(R_k)} \hat{H}^{\beta_k},
	\end{aligned}
	\end{align*}
	namely,
	\begin{align}\label{int_esti_iteration_2}
	\begin{aligned}
	||\hat{H}||_{L^{\beta_{k+1}}(B(R_{k+1}))} \leqslant \left(\frac{a_{10}\left(\beta_0^3+c_3\right)e^{c_1q_0}}{V^{2/n}}\right)^{1/\beta_k}\left( 16 \right)^{k/\beta_k}||\hat{H}||_{L^{\beta_{k}}(B(R_{k}))}.
	\end{aligned}
	\end{align}
	
	Then iterating (\ref{int_esti_iteration_2}) from $k = 1$ leads to 
	\begin{align}
	\begin{aligned}
	||\hat{H}||_{L^{\infty}(B(R/2))} &\leqslant 16^{\sum_{k=1}^{\infty} k/\beta_k}\left(\frac{a_{10}\left(\beta_0^3+c_3\right)e^{c_1q_0}}{V^{2/n}}\right)^{\sum_{k=1}^{\infty} 1/\beta_k} ||\hat{H}||_{L^{\beta_{1}}(B(3R/4))}\\
	&\leqslant e^{c_1}\left(16^{\frac{n^2}{4}}a_{10}^{\frac{n}{2}} \right)^{\frac{1}{\beta_1}}\left( \beta_0^3+c_3\right)^{\frac{n}{2\beta_1}}V^{-\frac{1}{\beta_1}}||\hat{H}||_{L^{\beta_{1}}(B(3R/4))}\\
	&\leqslant a_{13}V^{-\frac{1}{\beta_1}}||\hat{H}||_{L^{\beta_{1}}(B(3R/4))}\\
	&\leqslant C\frac{\left(1+\sqrt{K}R\right)^2}{R^2}.
	\end{aligned}
	\end{align}

	Finally, if $n=2$, Theorem \ref{sblv_thm} becomes
	\begin{align}\label{int_esti_sblv_n2}
	\begin{aligned}
	\left(\int_{\Omega} \hat{H}^{2(b+1)} \eta^{4}\right)^{\frac{1}{2}}\leqslant e^{C(1+\sqrt{K}R)}V^{-\frac{1}{2}}\left(R^2\int_{\Omega}\left|\nabla\left(\hat{H}^{b/2+1/2}\eta\right)\right|^2+ \int_{\Omega} \hat{H}^{b+1} \eta^{2}\right),
	\end{aligned}
	\end{align}
	where we take $n'=4$. Then By the same means, we have
	\begin{align*}
	||\hat{H}||_{L^{2\beta_0}(B_{3R/4})} \leqslant a_{14}V^{\frac{1}{2\beta_0}}\cdot\frac{q_0^2}{R^2},
	\end{align*}
	and
	\begin{align*}
	\begin{aligned}
	||\hat{H}||_{L^{\infty}(B(R/2))} \leqslant  a_{15}V^{-\frac{1}{2\beta_0}}||\hat{H}||_{L^{2\beta_0}(B(3R/4))}.
	\end{aligned}
	\end{align*}
	Thus, the desired estimate also holds for $n=2$.
\end{proof}

Then we give the direct applications of the gradient estimates.
\begin{corollary}[Harnack's inequality]\label{harnack}
	Under the same assumption in Theorem \ref{mainthm}, there exists a constant $C$ such that for any $x, y \in B(R)$,
	$$
	u(x) / u(y) \leq e^{C(1+\sqrt{K} R)} .
	$$
	If $K=0$, then we have a constant independent of $R$ such that
	$$
	\sup _{B(R)} u \leq C\inf _{B(R)} u .
	$$
	
\end{corollary}

\begin{proof} Under the same conditions in Theorem 1.1,
	let  $x$, $y\in B(R)$ be any two points with minimal geodesic $l$ connecting them.
	Then using the gradient estimate and the fact that $\operatorname{length}(l)\leqslant 2R$, we have
	\begin{align}
	\begin{aligned}
	\log u(x)- \log u(y) \leqslant\int_{l}|\nabla\log u|&\leqslant \int_{l}C\frac{\sqrt{K}R+1}{R}\\
	&\leqslant 2C(\sqrt{K}R+1).
	\end{aligned}
	\end{align}
	Therefore,
	$$
	u(x)  \leqslant e^{C(1+\sqrt{\kappa} R)}u(y) .
	$$
\end{proof}

\section{Nonlinear eigenvalue problems} 
 In this section, we continue to study nonlinear eigenvalue problems by the Cheng--Yau gradient estimate. Theorem \ref{nepthm} is directly deduced by the Liouville theorem as follows.
\begin{theorem}[Liouville theorem]\label{liou}
	Let $M$ be a complete and non-compact Riemannian manifold with non-negative Ricci curvature, and let $u$ be a bounded positive solution of (\ref{defharm}) with $\varphi$ and $\Psi$ satisfying the same assumption in Theorem \ref{mainthm}, then $\lambda = 0$. 
\end{theorem}
\begin{proof}
	When $K=0$ and $0< u\leqslant C$ is a bounded positive solution of 
	$$\Delta_\varphi(u) + \Psi(u^2,H)u = 0, $$ letting $R \to \infty$, we see $|\nabla u|=0.$  Consequently, $\Psi(u^2, 0) = 0 +  \frac{\lambda\psi(u^2)}{a(u^2)} = 0$. $u$ must be a positive constant and $\lambda = 0$.
\end{proof}

Then we will apply Theorem 1.1 to the nonlinear eigenvalue problem for a specific  $Q(u):=\Delta_{p}(\mathcal{P}(u))$, where $p>1$ and $\mathcal{P}(u) = \sum_{i=1}^{m}a_iu^{q_i}$ ($q_1<q_2<\cdots< q_m$) is a polynomial of $u$, satisfying $a_iq_i > 0$ for each $i=1,\cdots, m$ and
\begin{align}\label{p1pm}
	\frac{q_1^2(p-1)^2}{n-1} - \frac{(p-1)(q_m-q_1)^2}{2}>0
\end{align}

Now we cite a lemma in \cite{shen2025feasibilitynashmoseriterationchengyautype}.
\begin{lemma}[\cite{shen2025feasibilitynashmoseriterationchengyautype}]\label{lem}
	Let 
	$$f(t):= \frac{\sum_{i=1}^{m}a_ir_it^{r_i}}{\sum_{i=1}^{m}a_it^{r_i}}$$	
	where $a_i>0$ and $r_1<r_2<\cdots<r_m$, then for any $t\in[0,\infty)$, we have
	\begin{align}\label{lemp_1}
		r_{1}\leqslant f(t)\leqslant r_m
	\end{align}
\begin{align}\label{lemp_2}
	0\leqslant tf'(t) \leqslant (r_m-r_1)^2/2
\end{align}
\end{lemma}

The following theorem expend the one in \cite{he2023gradient} and \cite{shen2025feasibilitynashmoseriterationchengyautype} for more general operators.

\begin{theorem}
	For the nonlinear eigenvalue problems
	\begin{align}\label{nep_1}
		\Delta_{p}\left(\sum_{i=1}^{m}a_iu^{q_i}\right)+\lambda u^r = 0
	\end{align}
	on a complete and non-compact Riemannian manifold $M$ with Ricci curvature bounded from below by $\operatorname{Ric} \geqslant -K$ {where} $K\geqslant0$.
	We suppose that
	$p>1$, $q_1<q_2<\cdots< q_m$ satisfy (\ref{p1pm}). If
	\begin{align}\label{p1pm_1}
		\lambda\geqslant 0 \text{ and } r<\frac{(n+1)(p-1)q_1}{n-1} + \sqrt{\frac{4q^2_1(p-1)^2}{(n-1)^2} - \frac{2(q_m-q_1)^2(p-1)}{n-1}},
	\end{align}
	or
	\begin{align}\label{p1pm_2}
		\lambda\leqslant 0 \text{ and } r>\frac{(n+1)(p-1)q_m}{n-1} + \sqrt{\frac{4q^2_1(p-1)^2}{(n-1)^2} - \frac{2(q_m-q_1)^2(p-1)}{n-1}},
	\end{align}
	then the eigenfunction $u$ satisfies
	\begin{align}
	\frac{|\nabla u|}{u}\leqslant C(p,q_1,q_m,r)\frac{1+\sqrt{K}R}{R}
	\end{align}
	on $B(o,R)$.
	
	In particular, if $M$ is non-compact Riemannian manifold with non-negative Ricci curvature. For each non-zero $\lambda$, the eigenfunction $u$ is unbounded.
\end{theorem}

\begin{proof}
	It is easy to check that (\ref{nep_1}) is equivalent to the form of (\ref{defharm}) by choosing
	$$\varphi(t) = t^{(p-2)/2},$$ 
	$$a(s) = \left(\sum_{i=1}^{m}a_iq_i s^{(q_i-1)/2}\right)^{p-1},$$
	$$\psi(s) = \lambda s^{(r-1)/2}.$$
	Hence, we obtain that $\delta_{\varphi}\equiv p-2$, $\delta_{a} \equiv (q-1)(p-1)$, $\delta_{\psi} \equiv r-1$ and
	$$\delta_{a}(s)= \frac{(p-1)\sum_{i=1}^{m}a_iq_i(q_i-1)s^{(q_i-1)/2}}{\sum_{i=1}^{m}a_iq_i s^{(q_i-1)/2}}.$$
	Then, from Lemma \ref{lem}, the conditions (C1), (C2) follow that $p>1$ and 
	$$\gamma = \frac{4q^2_1(p-1)^2}{(n-1)^2} - \frac{2(q_m-q_1)^2(p-1)}{n-1}>0.$$
	
	Note that either $I = (0,+\infty)$ or $I = \emptyset$. The former case implies that 
	\begin{align}\label{pr1}
	\frac{n+1}{n-1}-\frac{r}{q_1(p-1)}\geqslant 0 \text{ when } \lambda\geqslant  0,
	\end{align}
	or 
	\begin{align}\label{pr2}
	\frac{n+1}{n-1}-\frac{r}{q_m(p-1)}\leqslant 0\text{ when } \lambda\leqslant  0.
	\end{align}
	The latter one holds if and only if
	$$\sup_{s\geqslant 0}\left(\frac{n+1}{n-1}\delta_a(s)-r\right)^2<\frac{4q^2_1(p-1)^2}{(n-1)^2} - \frac{2(q_m-q_1)^2(p-1)}{n-1},$$
	which infers that 
	\begin{align}\label{pr3_1}
	\frac{n+1}{n-1}q_m(p-1)-r<\sqrt{\frac{4q^2_1(p-1)^2}{(n-1)^2} - \frac{2(q_m-q_1)^2(p-1)}{n-1}},
	\end{align}
	and
	\begin{align}\label{pr3_2}
	\frac{n+1}{n-1}q_1(p-1)-r>-\sqrt{\frac{4q^2_1(p-1)^2}{(n-1)^2} - \frac{2(q_m-q_1)^2(p-1)}{n-1}}.
	\end{align}
	Combining (\ref{pr1}), (\ref{pr2}), (\ref{pr3_1}) and (\ref{pr3_2}), we obtain the desired results. 
\end{proof}
\begin{remark}
	In particular, if $\mathcal{P}(u) = u^q$, obiviously (\ref{p1pm}) holds. then (\ref{p1pm_1}) and (\ref{p1pm_2}) become
	$$\lambda\geqslant 0 \text{ and } r<\left(\frac{(n+1)q + 2|q|}{n-1}\right)(p-1)$$
	or
	$$\lambda\leqslant 0 \text{ and } r>\left(\frac{(n+1)q - 2|q|}{n-1}\right)(p-1).$$
	\begin{enumerate}
		\item Let $q=1$, then the result is same in \cite[Theorem 4.1]{he2023gradient}.
		%\item Let $p=2$ and $r=1$, our result not only expand the result in [?] for $q<0$, but also have larger admissible reigm for $q>0$. In [?], $q$ is required to satisfied
		%$$\left(\frac{n+1}{n-1}q-1\right)\leqslant\frac{2q^2}{(n-1)^2}.$$ 
		%As showed in (\ref{pr3}), it in fact holds for 
		%$$\left(\frac{n+1}{n-1}q-1\right)<\frac{4q^2}{(n-1)^2}.$$
		%Also their estimate is for weighted Riemannian manifolds, and we will show in our subsequent paper that this estimate is still true for general Finsler metric measure spaces.
		\item Let $p=2$ and $\lambda = r = 0$, our result also generalizes the Liouville theorems in \cite[Corollary 1.5 and Corollary 1.11]{huang2025gradient}, since we removed the requirement of positive lower bound of $u$.
	\end{enumerate}
\end{remark}

%\section*{Acknowledge}
%The authors are very grateful to the reviewers for their careful reviews and valuable comments. 

%\bibliography{refsab}
%\bibliographystyle{abbrv}

\end{document}